\documentclass[12pt]{amsart}%

\usepackage{amssymb}
\usepackage{amsmath}
\usepackage{amsfonts}
\usepackage{graphicx}%
\usepackage{verbatim}%
\newtheorem{theorem}{Theorem}[section]
\theoremstyle{plain}

\newtheorem{lemma}{Lemma}[section]

\newtheorem{proposition}{Proposition}[section]

\numberwithin{equation}{section}

\newcommand{\xionetwo}{\xi^1_{2}}
\newcommand{\calS}{\mathcal{S}}
\newcommand{\so}{\mathfrak{so}}
\newcommand{\qqand}{\qquad\text{and}\qquad}
\newcommand{\qand}{\quad\text{and}\quad}

\begin{document}
\title[Infinitesimal isometries along curves]{INFINITESIMAL\ ISOMETRIES\ ALONG\ CURVES\ AND\ GENERALIZED\ JACOBI\ EQUATIONS}
\author{R. L. Foote}
\address[Robert L. Foote]{ Department of Mathematics and Computer Science, Wabash College, Crawfordsville,
IN 47933, U.S.A.}
\email{footer@wabash.edu}
\author{C. K. Han$^{\ast}$}
\address[Chong-Kyu Han]{ Department of Mathematics, Seoul National University, Seoul
151-742, Korea}
\email{ckhan@snu.ac.kr}
\author{J. W. Oh$^{\ast}$}
\address[Jong-Won Oh]{ Department of Mathematics, Seoul National University, Seoul
151-742, Korea}
%\email[Jong-Won Oh]{jwoh@math.snu.ac.kr}
\thanks{$^\ast$ Partially supported by Korea Research Foundation grant C00007}
\thanks{To appear in J. Geometric Analysis.}
\date{4 September 2011}
\dedicatory{Jong-Won Oh passed away in May 2008. This paper is dedicated to his memory.}

\begin{abstract}
On a Riemannian manifold, a solution of the Killing equation is an infinitesimal isometry. Since the Killing equation is overdetermined, infinitesimal isometries do not exist in general. A completely determined prolongation of the Killing equation is a PDE on the bundle of 1-jets of vector fields. Restricted to a curve, this becomes an ODE that generalizes the Jacobi equation. A solution of this ODE is called an infinitesimal isometry along the curve, which we show to be an infinitesimal rigid variation of the curve. We define Killing transport to be the associated linear isometry between fibers of the bundle along the curve, and show that it is parallel translation for a connection on the bundle related to the Riemannian connection.  Restricting to dimension two, we study the holonomy of this connection, prove the Gauss-Bonnet theorem by means of Killing transport, and determine the criteria for local existence of infinitesimal isometries.
\end{abstract}

\maketitle

\section*{Introduction}

Let $(M^{n},g)$ be a smooth ($C^{\infty}$) manifold of dimension $n$ with
smooth Riemannian metric $g.$ An infinitesimal isometry, or a Killing field,
is a smooth vector field $X$ on $M$ whose flow $\phi_{t}:M\rightarrow M$ is an
isometry for each $t.$ Equivalently, $X$ is an infinitesimal isometry if and
only if it satisfies the Killing equation,
\begin{equation}
L_{X}g=0, \tag{1}\label{Killing equation}%
\end{equation}
where $L$ is the Lie derivative. Expressed in local coordinates or with respect to a local
frame, (\ref{Killing equation}) is a system of $n(n+1)/2$ linear partial differential
equations of first order for $n$ unknown functions. This is overdetermined if
$n\geq2,$ thus solutions do not exist in generic cases.

Our purpose is to introduce the notion of infinitesimal isometry along a
curve and Killing transport as a useful ODE analogs of (\ref{Killing equation}). The
situation is similar to that of parallel transport along a curve: the equation
$\nabla X=0$ for parallel vector fields on an open set of $M$ is an
overdetermined system of PDEs, generically with no solutions, which gives
rise to the first-order linear ODE $\nabla_{\gamma'(t)}X=0$ for a
vector field $X$ along a curve $\gamma$, and to the notion of parallel
transport along $\gamma.$ An infinitesimal isometry along $\gamma$ is given by
a second-order linear ODE that is the restriction to $\gamma$ of the completely determined
prolongation of (\ref{Killing equation}) to second order. In Section~\ref{prolongation} we discuss the notion
of a completely determined prolongation and its restriction to a curve in an abstract setting.

We present two different approaches to the second-order ODE that defines an
infinitesimal isometry along a curve. The first, presented in Section~\ref{inf isom}, is
the prolongation of (\ref{Killing equation}) and construction of a completely determined prolongation using the
moving frame method. The second approach, presented in Section~\ref{rigid var}, is from the
viewpoint of variation of curves. Given a curve $\gamma$, we consider a
smooth, rigid variation $\gamma_{\tau}$ of $\gamma=\gamma_{0}$. Intuitively,
one can think of a rigid variation as the motion in $M$ of a rigid, bent wire.
The variation involves both a translation, given by a vector field $X$ along
$\gamma$, and a rotation, given by a skew-symmetric $(1,1)$-tensor $A$ along
$\gamma$. The pair $(X,A)$ represents the 1-jet of an infinitesimal isometry along $\gamma$. We show that $X$ and $A$ satisfy
\begin{equation}
\nabla_{T}X=AT\qquad\text{and\qquad}\nabla_{T}A=R(T,X) \tag{2}\label{gen Jacobi}%
\end{equation}
along $\gamma$, where $T=\gamma'(t)$, which we think of as generalized Jacobi equations. (Indeed, if $\gamma$ is a geodesic, then $\gamma_{\tau}$ is a variation of geodesics, and (\ref{gen Jacobi}) easily implies that $X$ satisfies the classical Jacobi equation.) Conversely, we show that every solution of (\ref{gen Jacobi}) along $\gamma$ with $A$ skew-symmetric arises from such a rigid variation of $\gamma$. An inspection shows that (\ref{gen Jacobi}) is the restriction to $\gamma$ of the following second-order linear PDE system for $X$:
\begin{equation}\tag{3}\label{Jacobi prolongation}
\nabla X=A \text{ \ and \ } \nabla_Y A = R(Y,X),\ \ \text{or simply}\ \ \nabla_Y(\nabla X) =R(Y,X)
\end{equation}
for all vectors $Y$. When $\nabla X = A$ is skew, (\ref{Jacobi prolongation}) is a completely determined prolongation of (\ref{Killing equation}) to second order. Note that the equations in (\ref{gen Jacobi}) and (\ref{Jacobi prolongation}) do not require the skew-symmetry of $A$. We comment more on this in Sections~\ref{rigid var} and~\ref{connection}.

We define an \textsl{infinitesimal isometry along}
$\gamma$ to be a solution $(X,A)$ of (\ref{gen Jacobi}) when $A$ is skew. In addition, we define \textsl{Killing transport\/} to be the associated linear isometry between fibers of $TM \oplus \so(TM) \to M$ along $\gamma$, where $\so(TM) \to M$ denotes the bundle of skew endomorphisms of $TM$. Killing transport may be regarded as parallel
transport with respect to a connection $\tilde{\nabla}$ on
the bundle of skew 1-jets of vector fields, $TM\oplus\so(TM)\rightarrow M$. We investigate properties of this connection in Section~\ref{connection}. 

In Section~\ref{GBT} we use Killing transport to give another proof of the Gauss-Bonnet theorem.

In Section~\ref{existence} we discuss the local existence of solutions to (\ref{Killing equation}) in
dimension two. Our proof of the classical result in Theorem~\ref{class of inf isom} is based on discussions between the second author and R.~Bryant.

We would also like to thank the referee for many useful suggestions.

%%%%%%%%%%%%%%%%%%%%%%%%%%%%%%%%%%%%%%%%%%%%%%%%%%%%%%%%%%%%%%%%%%%%%%%

\section{Prolongation of overdetermined PDE systems\label{prolongation}}

Let $\Omega$ be an open subset of $\mathbb{R}^{n}.$ We consider a system of
linear partial differential equations of first order%
\begin{equation}
\sum_{\alpha=1}^{q}L_{\alpha}^{\lambda}u^{\alpha}=0,\quad\lambda=1,\dots,l
\label{OrigEqtn}%
\end{equation}
for unknown functions $u=(u^{1},\dots,u^{q})$ of independent variables $x$ $=$
$(x^{1},\dots,x^{n})$ $\in$ $\Omega,$ where each $L_{\alpha}^{\lambda}$ is a
linear partial differential operator of first order with coefficients that are smooth ($C^{\infty}$)
functions on $\Omega.$ We assume (\ref{OrigEqtn}) is overdetermined, that is,
$l>q.$ We discuss the case where, by differentiating (\ref{OrigEqtn}), one
obtains all the second-order partial derivatives of $u$ in terms of $u$ and
its first-order derivatives:%
\begin{equation}
u_{ij}^{\alpha}=H_{ij}^{\alpha}(x,u^{(1)}), \label{CompProlong}%
\end{equation}
for $\alpha=1,\dots,q$ and $i,j=1,\dots,n$. Here $u^{(1)}$ denotes $u$ and its
first-order partial derivatives and each $H_{ij}^{\alpha}$ is a linear
combination of $u^{\beta},$ $u_{k}^{\beta},$ for $\beta=1,\dots,q$ and
$k=1,\dots,n$, with $C^{\infty}$ coefficients. Equation (\ref{CompProlong}) is
called a \textsl{completely determined prolongation of (\ref{OrigEqtn}) to second order}. In this
case, we have that

%\bigskip

\begin{itemize}
\item[(i)] A solution $u$ is determined uniquely by its $1$-jet at a point,
and therefore the space of solutions is finite-dimensional.

\item[(ii)] If a solution $u$ is $C^{2}$, then $u$ is $C^{\infty}$
($C^{\omega}$ if the coefficients are $C^{\omega}$).

\item[(iii)] The conditions for solutions to exist can be found by checking
the Frobenius integrability conditions, or more generally, by the classical theory of Pfaffian systems.
\end{itemize}

%\bigskip

Properties (i) -- (iii) can be shown easily by defining a system of $1$-forms
on a subset of the first jet space: Assume that (\ref{OrigEqtn}) defines a
smooth submanifold $\calS$ of the first jet space $J^{1}(\Omega
,\mathbb{R}^{q}\mathbb{)=\{(}x,u^{(1)}\mathbb{)}:x\in\Omega\mathbb{\}}$. At
each point $(x,u^{(1)})\in\calS$ we consider a smooth distribution
$\mathcal{D}$\ of dimension $n$ annihilated by the differential $1$-forms%
\begin{equation}%
\begin{array}
[c]{l}%
\rule[-2mm]{0mm}{6mm}\theta^{\alpha}=du^{\alpha}-\sum_{i=1}^{n}u_{i}^{\alpha
}dx^{i}\\
\rule[-2mm]{0mm}{6mm}\theta_{i}^{\alpha}=du_{i}^{\alpha}-\sum_{j=1}^{n}%
H_{ij}^{\alpha}(x,u^{(1)})dx^{j}.
\end{array}
\label{EDSThetas}%
\end{equation}
Then the integral manifolds of the system
\begin{equation}
\left\{
\begin{array}
[c]{lc}%
\theta^{\alpha}=0 & \\
\theta_{i}^{\alpha}=0 & \\
dx^{1}\wedge\cdots\wedge dx^{n}\neq0 & \qquad\text{(independence condition)}%
\end{array}
\right.  \label{EDSsys}%
\end{equation}
are in one-to-one correspondence with the $C^{\infty}$ solutions of (\ref{OrigEqtn}), and (i) -- (iii) follow as
consequences. Our reference for these ideas is \cite{BCGGG}. A completely determined prolongation to third order has been studied in \cite{CH}.
To discuss the existence of solutions we use the following lemma, which is easy to prove (see \cite{H}).

\begin{lemma}
\label{MainLemma}Let $M$ be a smooth manifold of dimension $n$. Let
$\theta:=(\theta^{1},\dots,\theta^{s})$ be a set of independent 1-forms on $M$
and let $\mathcal{D}:= {<}\theta{>}^{\bot}$ be the $(n-s)$-dimensional
distribution annihilated by $\theta$. Suppose that $N$ is a submanifold of $M$
of dimension $n-r$, with $r\,{\leq}\, s,$ defined by $T_{1}=\dots=T_{r}=0$, where the
$T_{i}$ are smooth real-valued functions on $M$ such that $dT_{1}\wedge\cdots\wedge dT_{r}\neq0$ on $N$. Then the following are equivalent:

\begin{itemize}
\item[(i)] $\mathcal{D}$ is tangent to $N$.

\item[(ii)] $dT_{j}\equiv0$ $\operatorname{mod}$ $\theta$ at all points of $N$, for each
$j=1,\ldots,r.$

\item[(iii)] $i^{\ast}\theta= (i^{\ast}\theta^{1},\dots,i^{\ast}\theta^{s})$
has rank $s-r$, where $i:N\hookrightarrow M$ is the inclusion.
\end{itemize}
\end{lemma}

A key observation of this paper is that even though (\ref{EDSsys}) has no solutions in generic cases, given an initial condition, any smooth curve $\gamma:[a,b] \to \Omega$ has a unique lift $\tilde{\gamma}:[a,b] \to J^1(\Omega, \mathbb{R}^q)$ that is an integral curve  of (\ref{EDSsys}). The curve $\tilde\gamma$ is given by a solution of a system of linear ODE's, and the domain of $\tilde\gamma$ is all of $[a,b]$ since the ODEs are linear. If $\gamma$ is a closed curve, the quantity $\tilde{\gamma}(b)-\tilde{\gamma}(a)$ measures the total torsion (non-integrability) of
(\ref{EDSsys}) along $\gamma$.

\bigskip\noindent{\bf The Killing Equation.} 
Now consider the situation for the Killing equation (\ref{Killing equation}). Let $x=(x^{1},\dots,x^{n})$
be a local coordinate system of $M^{n}$ and let $g_{ij}=g\left(
\frac{\partial}{\partial x^{i}},\frac{\partial}{\partial x^{j}}\right)  $. Let
$X=\sum u^{\alpha}\frac{\partial}{\partial x^{\alpha}}$. Then in terms of
these coordinates, (\ref{Killing equation}) can be written as
\begin{equation}
\sum_{\alpha=1}^{n}\left(  g_{\alpha j}\frac{\partial u^{\alpha}}{\partial
x^{i}}+g_{i\alpha}\frac{\partial u^{\alpha}}{\partial x^{j}}+\frac{\partial
g_{ij}}{\partial x^{\alpha}}u^{\alpha}\right)  =0, \label{Killing loc coords}%
\end{equation}
for each $i,j=1,\dots,n.$ Since $g_{ij}=g_{ji}$, the number of equations in
(\ref{Killing loc coords}) is $n(n+1)/2$, and therefore (\ref{Killing loc coords}) is overdetermined when $n\geq2.$ A completely determined prolongation to second order can be obtained by differentiating
(\ref{Killing loc coords}) with respect to each coordinate and by solving linear algebraic
equations for all the second-order partial derivatives of $u$. In this paper,
we present two different approaches to coordinate-free calculations of this prolongation: One is by using the moving frame method in Section~\ref{inf isom};
the other is by using rigid variations and the Levi-Civita covariant
derivative operator $\nabla$ in Section~\ref{rigid var}.

%%%%%%%%%%%%%%%%%%%%%%%%%%%%%%%%%%%%%%%%%%%%%%%%%%%%%%%%%%%%%%%%%%%%%%%%%%%%%%%%%%%

\section{Infinitesimal isometries on Riemannian manifolds\label{inf isom}}

In this section we discuss the completely determined prolongation of (\ref{Killing equation}) to second order by the method of moving frames. Let $e_{i},$ $i=1,\dots,n,$ be a local, orthonormal
frame on a Riemannian manifold $(M^{n},g)$ and let $\omega^{i},$ $i=1,\dots,n,$ be the dual coframe. Then $g=\sum_{i=1}^{n}\omega^{i}\circ\omega^{i},$ where $\phi\circ\eta$ $:=$ $\frac{1}{2}(\phi\otimes\eta$
$+$ $\eta\otimes\phi)$ is the symmetric product of $1$-forms. Recall also that
there exist uniquely determined $1$-forms $\omega^{i}_{j},$ $i,j=1,\dots,n$ (Levi-Civita connection) and the curvature tensor $R_{ijkl}$
satisfying%
\begin{align}
d\omega^{i}  &  =\omega^{k}\wedge\omega_{k}^{i},\label{connection forms}\\
d\omega^{i}_{j} & =\omega_{j}^{k}\wedge\omega_{k}^{i}+\frac{1}{2}R_{ijkl}\,\omega^{k}\wedge\omega^{l} 
   \label{curvature tensor}
\end{align}
with the symmetries
$$
\begin{aligned}
\omega^i_j + \omega^j_i &= 0 \\ \nonumber
R_{sijk} + R_{sjki} + R_{skij} &= 0  \qquad \text{(1st Bianchi identity)}.
\end{aligned}
$$
(Since we are working with an orthonormal frame, we may use all lowered indices for the curvature tensor.)

\begin{proposition}
\label{MainSystem} Equation (\ref{Killing equation}) for an infinitesimal isometry $X=\xi^{j}e_{j}$ on $M$ admits a completely determined prolongation to second-order as follows:%
\begin{equation}
\begin{aligned}
\theta^{i}  &  =d\xi^{i} + \xi^{k}\,\omega_{k}^{i} - \xi_{k}^{i}\,\omega^{k}
\text{\quad with\quad} \xi_{k}^{j}+\xi_{j}^{k}=0,\\
\theta_{j}^{i} & = d\xi_{j}^{i} - \xi_{k}^{i}\,\omega_{j}^{k} + \xi_{j}^{k}\,\omega_{k}^{i} + \xi^{k}R_{ijkl}\,\omega^{l},%\nonumber
\end{aligned}
\label{ProlongFrame}
\end{equation}
where the system $\{\theta^{i},\theta_{j}^{i}\}$ is defined on the
bundle $J^{1}(TM)\rightarrow M.$
\end{proposition}

\begin{proof}
Suppose $X$ is an infinitesimal isometry. Then
\begin{equation}
0=L_{X}g=2\sum_{i=1}^{n}\omega^{i}\circ L_{X}\omega^{i}. \label{Killing frame}%
\end{equation}
We have the identity
\begin{align}
L_{X}\omega^{i}  &  =d(X\lrcorner\,\omega^{i}) + X\lrcorner\, d\omega
^{i}\nonumber\\
&  =d\xi^{i}+\omega^{k}(X) \omega_{k}^{i}-\omega_{k}^{i}( X)
\omega^{k}\label{Lie der ident}\\
&  =d\xi^{i}+\xi^{k}\omega_{k}^{i}-\omega_{k}^{i}(X)
\omega^{k}\nonumber\\
&  =\xi_{k}^{i}\omega^{k}-\omega_{k}^{i}( X)\omega
^{k},\nonumber
\end{align}
where we define the quantity $\xi_{k}^{i}$ by
\begin{equation}
d\xi^{i}+\xi^{k}\omega_{k}^{i}=\xi_{\ k}^{i}\omega^{k}.
\label{xi ik def}%
\end{equation}
Substituting (\ref{Lie der ident}) for $L_{X}\omega^{i}$ in (\ref{Killing frame}),
\begin{align*}
0  &  =\sum_{i,k=1}^{n}\omega^{i}\circ(\xi_{k}^{i}\omega^{k}-\omega_{k}^{i}(X)  \omega^{k})\\
&  =\sum_{i}\xi_{i}^{i}\omega^{i}\circ\omega^{i}+\sum_{j<k}%
(\xi_{k}^{j}+\xi_{j}^{k})\omega^{j}\circ\omega^{k}
+\sum_{j<k}(\omega_{j}^{k}(X)+\omega_{k}^{j}(X))\omega^{j}\circ\omega^{k}.
\end{align*}
Since $\omega_{j}^{k}+\omega_{k}^{j}=0,$ we have%
\[
\xi_{k}^{j}+\xi_{j}^{k}=0.
\]
Differentiate $\ \ d\xi^{i} = -\xi^{k}\,\omega_{k}^{i} + \xi_{k}^{i}\,\omega^{k}$ \ to obtain%
\begin{align*}
0  &  =-d\xi^{k}\wedge\omega_{k}^{i}-\xi^{k}d\omega_{k}^{i}+d\xi_{k}^{i}\wedge\omega^{k}+\xi_{k}^{i}\,d\omega^{k}\\
&  =(\xi^{l}\omega_{l}^{k}-\xi_{l}^{k}\omega^{l})\wedge\omega_{k}^{i}
    - \xi^{k}(\omega_{k}^{l}\wedge\omega_{l}^{i} + \frac{1}{2}R_{ikpq}\,\omega^{p}\wedge\omega^{q})\\
&  \qquad\qquad + d\xi_{k}^{i}\wedge\omega^{k} + \xi_{k}^{i}\,\omega^{l}\wedge\omega_{l}^{k}\\
&  =(d\xi_{k}^{i} + \xi_{k}^{l}\,\omega_{l}^{i} - \xi_{l}^{i}\,\omega_{k}^{l} 
     - \frac{1}{2}\xi^{l}R_{iljk}\,\omega^{j})\wedge\omega^{k}.
\end{align*}
\newline Then by Cartan's lemma,%
\begin{equation}\label{def Cijk}
d\xi_{k}^{i} + \xi_{k}^{l}\,\omega_{l}^{i} - \xi_{l}^{i}\,\omega_{k}^{l}
  - \frac{1}{2}\xi^{l}R_{iljk}\,\omega^{j} = C_{kj}^{i}\omega^{j} 
\end{equation}
with $C_{ks}^{i} = C_{sk}^{i}.$ Switching $i$
and $k$ we have%
\begin{equation}\label{switch ik}
d\xi_{i}^{k} + \xi_{i}^{l}\,\omega_{l}^{k} - \xi_{l}^{k}\,\omega_{i}^{l}
  - \frac{1}{2}\xi^{l}R_{klji}\,\omega^{j} = C_{ij}^{k}\omega^{j} 
\end{equation}
Since $\omega_{k}^{j}=-\omega_{j}^{k}$ and $\xi_{%
k}^{j}=-\xi_{j}^{k},$ the sum of (\ref{def Cijk}) and (\ref{switch ik}) gives%
\[
-\frac{1}{2}\xi^{l}(R_{iljk} + R_{klji})\,\omega^{j}=(C_{ks}^{i}+C_{is}^{k})\,\omega^{s}.%
\]
Rearranging indices using $C_{ks}^{i}=C_{sk}^{i}$ we have%
\[
-\frac{1}{2}\xi^{s}(R_{isjk} - R_{ksij}) = C_{jk}^{i} + C_{ij}^{k}.
\]
Combine the equations with $i,j,k$ permuted to obtain%
\begin{equation}
-\frac{1}{2}\xi^{s}(R_{jski} - R_{ksij}) = C_{jk}^{i}. \label{combined}%
\end{equation}
By substituting (\ref{combined}) for $C_{ij}^{k}$ into (\ref{switch ik}) and using $C_{jk}^{i}=C_{kj}^{i}$, we have
\begin{align*}
d\xi_{k}^{i} + \xi_{k}^{l}\,\omega_{l}^{i} - \xi_{l}^{i}\,\omega_{k}^{l}  
&  = \frac{1}{2}\xi^{l}R_{iljk}\,\omega^{j} - \frac{1}{2}\xi^{s}(R_{jski} - R_{ksij})\,\omega^{j}\\
&  = \frac{1}{2}\xi^{s}\left( - R_{sijk} + R_{sjki} - R_{skij}\right)  \omega^{j}\\
&  = \xi^{s}R_{sjki}\omega^{j}\qquad\text{(by the 1st Bianchi identity)}\\
&  = -\xi^{s}R_{iksj}\omega^{j}.
\end{align*}
Thus we define%
\[
\theta_{j}^{i} := d\xi_{j}^{i} + \xi_{j}^{k}\,\omega_{k}^{i} - \xi_{k}^{i}\,\omega_{j}^{k} + \xi^{k}R_{ijkl}\,\omega^{l}.
\]
Then%
\[
\theta^{i}=0,\qquad\theta_{j}^{i}=0,\qquad i,j=1,\dots,n
\]
is a completely determined prolongation to second order for $X=\xi^{j}e_{j}$.
\end{proof}

In the next section, we take another viewpoint of the prolongation of (1).

\section{Rigid Variations of a Curve and Generalized Jacobi Equations\label{rigid var}}

Let $M^{n}$ be a Riemannian manifold. Let $\gamma:[0,L]\rightarrow
M$ be a $C^{2}$ curve parameterized by arc length. In this section we define a
rigid variation of $\gamma$, show that a rigid variation gives rise to a type of
generalized Jacobi field along $\gamma$, and that every suitable generalized
Jacobi field along $\gamma$ arises from such a variation. For technical simplicity we assume that $M$ is complete, although the results in this section depend only on the geometry of $M$ in a neighborhood of $\gamma$.

In order to specify the rigidity of $\gamma$, we choose a relatively parallel
frame along $\gamma$ following Bishop \cite{B}. Let $T=\gamma'(t)$ and
extend $T_{\gamma(0)}$ to an orthonormal frame $T_{\gamma(0)},$ $N_{2}%
,\dots,N_{n}$ of $T_{\gamma(0)}M.$ Extend each $N_{i}$ along $\gamma$ by
parallel translation in the normal bundle along $\gamma$ by the connection induced from $\nabla.$ Then $T,$ $N_{2},\dots,N_{n}$ remain orthonormal and
$\nabla_{T}N_{i}$ is a multiple of $T.$ We shall call this frame a
\textsl{Bishop frame} along $\gamma.$

Define geodesic curvature functions $\kappa_{2},\dots,\kappa_{n}$ 
by $\nabla_{T}T = \sum_{2}^{n} \kappa_{i}N_{i}$. It is easily
seen that $\nabla_{T}N_{i} = -\kappa_{i}T$. In fact, the initial frame
at $\gamma(0)$ and the geodesic curvature functions uniquely determine
$\gamma$, as described by the following theorem.

\begin{theorem}
\label{BishopFrameTheorem} Let $M^{n}$ be a complete, Riemannian manifold. Let $T,$ $N_{2},\dots,N_{n}$ be an orthonormal
frame for $T_{p}M,$ and let $\kappa_{2},$ $\dots,$ $\kappa_{n}$
$:[0,L]$ $\rightarrow$ $\mathbb{R}$ be continuous. Then there exists a unique
$C^{2}$ curve $\gamma:[0,L]$ $\rightarrow$ $M$, parameterized by arc length,
and unique extensions of $T,$ $N_{2},\dots,N_{n}~$\ to a Bishop frame along
$\gamma$ such that $\gamma(0)=p,$ $\gamma'(t)=T,$ and $\nabla_{T}T$
$=$ $\sum_{2}^{n}$ $\kappa_{i}N_{i}.$
\end{theorem}

\noindent The proof is by existence and uniqueness for the ODE system $\dot
{T}=\sum_{2}^{n}$ $\kappa_{i}N_{i}$ and $\dot{N}_{i}=-\kappa_{i}T$ on
the orthonormal frame bundle $\mathcal{O(}TM)$ $\rightarrow$ $M.$ See \cite{B}
and \cite[pg. 121]{O} for proofs in $\mathbb{R}^{n}$, which are easily adapted
to the case needed here.

Given the curve $\gamma$ in $M,$ a \textsl{rigid variation} of $\gamma$ is a
one-parameter family of curves $\{\gamma_{\tau}\},$ $|\tau|<\epsilon,$ with
$\gamma=\gamma_{0},$ such that each curve has its own Bishop frame and all the
curves in the family have the same geodesic curvature functions. Specifically,
let $c$ $:$ $(-\epsilon,\epsilon)$ $\rightarrow$ $M$ be a $C^{1}$ curve such
that $c(0)$ $=$ $\gamma(0).$ Let $T$ $=$ $\gamma'(0),N_{2},\dots
,N_{n}$ be an orthonormal frame in $T_{\gamma(0)}M.$ Extend this frame to an orthonormal frame along $c$ in an arbitrary $C^{1}$ manner, and to a Bishop frame along $\gamma.$
Define the functions $\kappa_{2},\dots,\kappa_{n}$ as above. Then by
Theorem~\ref{BishopFrameTheorem}, for each $\tau$ there exists a unique curve
$\gamma_{\tau}$ and a Bishop frame $T_{\tau},$ $(N_{2})_{\tau},\dots,$
$(N_{n})_{\tau}$ along $\gamma_{\tau}$ such that $\gamma_{\tau}(0)=c(\tau),$
$\gamma_{\tau}'(t)=T_{\tau}$ and $\nabla_{T_{\tau}}T_{\tau}$ $=$
$\sum_{2}^{n}$ $\kappa_{i}(N_{i})_{\tau}.$ The variation field of
$\gamma_{\tau}$ is the vector field $X$ along $\gamma$ defined by
$X_{\gamma(t)}=\left.  \frac{\partial}{\partial\tau}\right\vert _{0}%
\gamma_{\tau}(t).$

\begin{theorem}
\label{variation thm} Let $\gamma:[0,L]\to M$ be a $C^2$ curve parameterized by arc length. Let $\{\gamma_{\tau}\}$ be a rigid variation of
$\gamma$, and let $X$ be its variation field. Then there is a skew-symmetric
$(1,1)$ tensor $A$ along $\gamma$ such that $X$ and $A$ satisfy the
generalized Jacobi equations%
\begin{equation}
\nabla_{T}X=AT\qquad\text{and\qquad}\nabla_{T}A=R(T,X). \label{GenJacEqtn}%
\end{equation}
Conversely, every solution of (\ref{GenJacEqtn}) with $A$ skew-symmetric
arises from a rigid variation of $\gamma.$
\end{theorem}

\noindent{\bf Definition.} A solution $(X,A)$ of (\ref{GenJacEqtn}) with $A$ skew is called an \textsl{infinitesimal isometry along\/} $\gamma$. For $s,t\in[0,L]$, the linear isomorphism $T_{\gamma(s)}M \oplus\so(T_{\gamma(s)}M) \to T_{\gamma(t)}M \oplus\so(T_{\gamma(s)}M)$ given by $(X(s),A(s)) \mapsto (X(t),A(t))$ is called \textsl{Killing transport along\/} $\gamma$, where $\so(T_pM)$ denotes the skew endomorphisms of $T_pM$.

\smallskip
Note that if $A$ is skew-symmetric at one point of $\gamma$, it will be skew-symmetric
all along $\gamma$ since $R(T,X)$ is skew-symmetric. Kostant derives the equations (\ref{GenJacEqtn}) in \cite[pg.~535]{Kostant}. He notes that a vector field $X$ is an infinitesimal isometry (Killing field) if and only if $X$ and $A = \nabla X$ satisfy (\ref{GenJacEqtn}) along all differentiable curves.

\begin{proof}
[Proof of Theorem \ref{variation thm}]The proof is similar to the development
of the Jacobi equation in \cite[pg. 14]{CE}. Define $\Gamma:(-\epsilon
,\epsilon)\times\lbrack0,L]\rightarrow M$ by $\Gamma(\tau,t)=\gamma_{\tau
}(t).$ The vector fields $T_{\tau}$ and $(N_{i})_{\tau}$ along the curves
$\gamma_{\tau}$ form vector fields $T=\Gamma_{\ast}(\partial/\partial t)$ and $N_{i}$ along $\Gamma.$ The variation field $X$ extends to a vector field along
$\Gamma$ by setting $X_{(\tau,t)}$ $=$ $\Gamma_{\ast}\left(  \partial
/\partial\tau\right)  ,$ that is, $X$ $=$ $\frac{\partial}{\partial\tau}%
\gamma_{\tau}(t).$ Note that $\nabla_{T}X=\nabla_{X}T$ since $\nabla
_{T}X-\nabla_{X}T=\Gamma_{\ast}\left(  [\partial/\partial t,\partial
/\partial\tau]\right)  =0.$

Define the $(1,1)$ tensor $A$ by $AT=\nabla_{X}T$ and $AN_{i}=\nabla_{X}%
N_{i},$ that is, $A$ is the transverse derivative of the Bishop frame. Note that $A$ is skew-symmetric since the
frame is orthonormal. We have $\nabla_{T}X$ $=$ $\nabla_{X}T$ $=$ $AT,$ the
first equation of (\ref{GenJacEqtn}). 

To prove that $(\nabla_T A)T = R(T,X)T$, we have
$$
\begin{aligned}
(\nabla_T A)T &= \nabla_T(AT) - A(\nabla_T T) 
\textstyle= \nabla_T(\nabla_X T) - A\left(\sum \kappa_i N_i\right) \\
&\textstyle= R(T,X)T + \nabla_X\nabla_T T - \sum \kappa_i AN_i \\
&\textstyle= R(T,X)T + \nabla_X\left(\sum \kappa_i N_i\right) - \sum \kappa_i \nabla_X N_i \\
&\textstyle= R(T,X)T + \sum (\nabla_X\kappa_i) N_i.
\end{aligned}
$$
The geodesic curvature functions $\kappa_{i}$ do not depend on $\tau$, and so $\nabla_X\kappa_i = 0$ and $(\nabla_T A)T = R(T,X)T$ follows. The proof that $(\nabla_{T}A)N_{i}=R(T,X)N_{i}$ is similar. Thus we have $\nabla_{T}A=R(T,X),$ which proves the second equation in (\ref{GenJacEqtn}).

To prove the converse, suppose $X$ and $A$ satisfy (\ref{GenJacEqtn}) along
$\gamma,$ where $A$ is skew-symmetric. Extend $T_{\gamma(0)}$
to an orthonormal frame $T,N_{2},\dots,N_{n}$ of $T_{\gamma(0)}M.$ Let $c$ be
a curve in $M$ such that $c'(0)=X_{\gamma(0)}.$ Extend $T$ and $N_{i}$
to an orthonormal frame along $c$ that satisfies $\nabla_{X}T=AT$ and
$\nabla_{X}N_{i}=AN_{i}$ at $\gamma(0),$ which requires the skew-symmetry of
$A$. The process above then defines a rigid variation of $\gamma$, which
yields a solution $(\tilde{X},\tilde{A})$ of (\ref{GenJacEqtn}) along $\gamma$
with the same initial conditions as $(X,A).$ Thus $(\tilde{X},\tilde
{A})=(X,A)$ by the uniqueness theorem for ODEs, and so $(X,A)$ arises from a
rigid variation of $\gamma.$
\end{proof}

\noindent\textbf{Remarks.}
To justify calling (\ref{GenJacEqtn}) generalized Jacobi equations, suppose that $X$ and $A$ satisfy (\ref{GenJacEqtn}) and that $\gamma$ is a geodesic. Then $\nabla_T T = 0$ and we have
$$
\nabla_T^2 X = \nabla_T(AT) = (\nabla_T A)T + A(\nabla_T T) = R(T,X)T,
$$
and so $X$ satisfies the classical Jacobi equation.

Expressed in terms of a local orthonormal frame $e_1$, \dots, $e_n$ and its dual frame $\omega^1$, \dots, $\omega^n$, it is easily seen that a solution of (\ref{GenJacEqtn}) is a solution of the system (\ref{ProlongFrame}) restricted to $\gamma$:
\begin{equation}
\begin{aligned}
d\xi^{i}/dt &= -\xi^k\omega_k^i(\gamma'(t)) + \xi_k^i\omega^k(\gamma'(t)) \\
d\xi_j^i/dt 
&= \xi_k^i\omega_j^k(\gamma'(t)) - \xi_j^k\omega_k^i(\gamma'(t)) - \xi^k R_{ijks}\omega^s(\gamma'(t)), 
\end{aligned}
\label{GenJacEqtnFrame}
\end{equation}
in which $X = \xi^k e_k$ and the coefficients of $A$ are $\{\xi_{j}^{k}\}$.

Equations (\ref{GenJacEqtn}) are the restrictions to $\gamma$ of the equation
\begin{equation}
\nabla_Y(\nabla X) =R(Y,X) \label{ProlongCoordFree}%
\end{equation}
for all vectors $Y$ (c.f.~\cite[pg.~235]{KN}). Expressed in terms of the local frame, (\ref{ProlongCoordFree}) becomes
$\theta^{i} = 0$, $\theta^{i}_{j} = 0$, where $\theta^{i}$ and $\theta^{i}_{j}$ are given in (\ref{ProlongFrame}). Thus (\ref{ProlongCoordFree})
is another expression of the completely determined prolongation to second order for the
Killing equation (\ref{Killing equation}). Note that if $X$ is a global infinitesimal isometry, then it is a Jacobi field along every geodesic, which readily follows from (\ref{ProlongCoordFree}).

As the proof of the theorem shows, the tensor $A$ describes the infinitesimal rotation of the frame, and so it is a measure of the rotation of the rigid variation. If $X$ extends to an infinitesimal isometry on a neighborhood of $\gamma$, the flow of $X$ preserves the bundle of orthonormal frames, and $A$ represents the derivative of the flow on this bundle. This easily implies that $A = \nabla X$, which justifies thinking of $(X,A)$ as the 1-jet of an infinitesimal isometry on $\gamma$.

Note that nothing in equations (\ref{GenJacEqtn}) requires that $A$ be skew-sym\-met\-ric.
Relaxing this requirement (equivalently, relaxing the condition $\xi_{k}^{j}+\xi_{j}^{k}=0$ in (\ref{ProlongFrame})) results in the notion
of an infinitesimal affine transformation along $\gamma$ (recall that an
infinitesimal affine transformation on a manifold with connection is a vector
field whose flow preserves the connection \cite[pg.~230]{KN}).

%%%%%%%%%%%%%%%%%%%%%%%%%%%%%%%%%%%%%%%%%%%%%%%%%%%%%%%%%%%%%%%%%%%%%%%%%%%%%%%%%%%

\section{Killing transport as parallel transport \\ for a connection on $J^1(TM)\to M$\label{connection}}

Let $\gamma:[0,L]\rightarrow M$ be a a $C^2$ curve parameterized by arc length. A solution $(X,A)$ of (\ref{GenJacEqtn}) along $\gamma$ (with no assumption that $A$ is skew) may be regarded as parallel transport with respect to a connection on the bundle of 1-jets of vector fields $J^{1}(TM)=TM\oplus\text{End}(TM)\rightarrow M$ related to the Levi-Civita connection. Indeed, if $\tilde{\nabla}$ is defined by
\[
\tilde{\nabla}_{Y}(X,A)=\nabla_{Y}(X,A)-(AY,R(Y,X))=(\nabla_{Y}X-AY,\nabla
_{Y}A-R(Y,X))
\]
for sections $(X,A)$ of $J^{1}(TM)$, it is easy to show that $\tilde{\nabla}$
is a covariant derivative (or Koszul connection \cite{Spivak}). 

By Theorem~\ref{variation thm}, the parallel transport for $\tilde{\nabla}$ preserves the sub-bundle $TM\oplus\so(TM)\rightarrow M$ (on which it is called Killing transport), and so $\tilde{\nabla}$ restricts to a connection on this sub-bundle.  We note that this connection on this sub-bundle and its relation to Killing fields have been independently studied by other authors recently \cite{Atkins,ConsoleOlmos,Eastwood1,Eastwood2}.  In particular, they note that parallel sections of the sub-bundle correspond to Killing fields on $M$. It is not difficult to see that parallel sections of the larger bundle $J^{1}(TM)\rightarrow M$ correspond to infinitesimal affine transformations on $M$. 

In this section we study the properties of $\tilde{\nabla}$ on this sub-bundle in dimension two when $M$ is oriented, in which case $TM\oplus\so(TM) \cong TM\oplus\mathbb{R}$. Let $T = \gamma'$. Let $N$ be the unit normal vector field along $\gamma$ such that $T$, $N$ form an oriented frame. Since $A$ is skew-symmetric, it is a scalar multiple of $J$, where $JT=N$ and $JN=-T$. In fact, it is easily seen that $A = -\xionetwo J$, where $\xi_{i}^{j}$ is defined in (\ref{xi ik def}) for an arbitrary, oriented, orthonormal frame $e_1$, $e_2$. It follows that $\xionetwo$ is independent of the choice of oriented frame in (\ref{xi ik def}). Then equations (\ref{GenJacEqtn}) become%
\begin{equation}
\nabla_{T}X = -\xionetwo N\quad\text{and\quad}
d\xionetwo/dt = K\left\langle N,X\right\rangle = K\,d\alpha(T,X), \label{GenJacEqtn2D}%
\end{equation}
where $K = \left\langle R(T,N)N,T\right\rangle = R_{1212}$ is the Gaussian curvature of $M$ and $d\alpha$ is the area form. Relative to an oriented, orthonormal frame $e_1$, $e_2$ we have, either from (\ref{GenJacEqtn2D}) or (\ref{GenJacEqtnFrame}), that
\begin{equation}%
\begin{aligned}
d\xi^{1}/dt &= -\xi^{2}\omega_{2}^{1}(\gamma') + \xionetwo\omega^{2}(\gamma')\\
d\xi^{2}/dt &= \xi^{1}\omega_{2}^{1}(\gamma') - \xionetwo\omega^{1}(\gamma')\\
d\xionetwo/dt &= - (K{\circ}\gamma)\,\xi^{1}\omega^{2}(\gamma') + (K{\circ}\gamma)\,\xi^{2}\omega^{1}(\gamma'),
\end{aligned}
\label{SysOnCurve2D}%
\end{equation}
where $X = \xi^{1}e_1 + \xi^{2}e_2$. The system takes a particularly nice form when $e_1 = T$ and $e_2 = N$:
\begin{equation}
\begin{matrix}
{d\xi_T}/{dt} & = & & \kappa\xi_N \\
{d\xi_N}/{dt} & = & -\kappa\xi_T & & & -\xionetwo  \\
{d\xionetwo}/{dt}     & = & & (K{\circ}\gamma)\xi_N, \\
\end{matrix} \label{SysOnCurve2DSpecialFrame}
\end{equation}
where $X = \xi_{T}T + \xi_{N}N$ and $\kappa$ is the geodesic curvature of $\gamma$.

Let $Q(t)$ be the matrix for the system (\ref{SysOnCurve2D}), that is,
$$
Q(t)=%
\begin{pmatrix}
0 & -\omega_{2}^{1}(\gamma'(t)) & \omega^{2}(\gamma'(t)) \\
\omega_{2}^{1}(\gamma'(t)) & 0 & -\omega^{1}(\gamma'(t)) \\
-K(\gamma(t))\,\omega^{2}(\gamma'(t)) & K(\gamma(t))\,\omega^{1}(\gamma'(t)) & 0
\end{pmatrix}.
$$
If $U:[0,L]\to GL(\mathbb R^3)$ solves 
\begin{equation}
U'(t) = Q(t)U(t) \qqand U(0)=I, \label{matrix solution}
\end{equation}
then
$$
\big(  \xi^1(t), \xi^2(t), \xionetwo(t)\big)^T 
= U(t)\,\big(  \xi^1(0), \xi^2(0), \xionetwo(0)\big)^T
$$
is the solution to (\ref{SysOnCurve2D}). Since $\mathop{\mathrm{tr}}Q(t) = 0$, then $\mathop{\mathrm{det}}U(t) = 1$, and so $U(t) \in SL(\mathbb R^3)$. 
If $K$ is not constant on $\gamma$, a simple computation shows that the smallest Lie algebra containing every $Q(t)$ is $\mathfrak{sl}(\mathbb R^3)$. It follows that $SL(\mathbb R^3)$ is the smallest group containing every $U(t)$ in the general case.

\bigskip\noindent{\bf Curvature.} 
In two dimensions the curvature tensor of $\tilde{\nabla}$ can be shown to be
\begin{equation}
\tilde{R}(Y,Z)(X,A)=(0,-dK(X)\,d\alpha(Y,Z)J), \label{curvature}%
\end{equation}
where $d\alpha$ is the area form on $M$. Thus the curvature of $\tilde{\nabla}$ is essentially $dK$. It follows that non-constant Gaussian curvature is the main obstruction to the existence of local infinitesimal isometries. This will play a role in Section~\ref{existence}.

More generally, for the connection $\tilde{\nabla}$ on the full bundle $J^1(TM) = TM \oplus \text{End}(TM)\to M$ with $M$ of arbitrary dimension, the curvature tensor is
$$
\tilde R(Y,Z)(X,A) = (0, (\nabla_X R)(Y,Z) + [R(Y,Z),A] + R(Y,AZ) + R(AY,Z)),
$$
where $R$ is the curvature tensor of $\nabla$. (The derivation is a straightforward exercise using both Bianchi identities.) This curvature tensor is the main obstruction to the existence of local infinitesimal affine transformations.

\bigskip\noindent{\bf Holonomy.} 
Around a closed curve, the holonomy for Killing transport measures the
non-integrability of (\ref{Killing equation}). This is analogous to the familiar fact that the holonomy for ordinary parallel transport around a closed curve measures the non-integrability of the equation $\nabla X=0$.
For a loop $\gamma$ in $M$ with $\gamma(0)=\gamma(L)$, the holonomy around $\gamma$ is $U(L)$, where $U$ is defined by (\ref{matrix solution}). The holonomy is trivial if $U(L) = I$. The holonomy is trivial for the initial condition $\xi_0 = (\xi^1(0), \xi^2(0), \xionetwo(0))$ if $U(L)\xi_0^{\ T} = \xi_0^{\ T}$.

For example, suppose that $M \subset \mathbb R^3$ is a surface of revolution about some line $\ell$. If $\theta$ measures the angle of rotation around $\ell$, then $\partial/\partial\theta$ restricts to an infinitesimal isometry on $M$. Let the circle $\gamma$ be an integral curve of $\partial/\partial\theta$ on $M$. Since the geodesic curvature $\kappa$ of $\gamma$ is constant, $\xi_T$ can be eliminated in (\ref{SysOnCurve2DSpecialFrame}), yielding
$$
{d^2\xi_N}/{dt^2} = -(\kappa^2 + K)\xi_N  \quad\text{and}\quad {d\xionetwo}/{dt} = K\xi_N.
$$
Note that if $c$ is constant, then
$$
\xi_T = c, \qquad \xi_N = 0, \qquad \xionetwo = -\kappa c
$$
is a solution. Thus the initial condition $(c,0,-\kappa c)$ has trivial holonomy, which is expected since $M$ admits an infinitesimal isometry. The value $c = L/(2\pi)$ corresponds to the infinitesimal isometry $\partial/\partial\theta$.

An initial condition that is not a multiple of $(1,0,-\kappa)$ will lead to a non-constant solution of (\ref{SysOnCurve2DSpecialFrame}). Since $K$ is constant on $\gamma$, the solution has period $2\pi/\sqrt{\kappa^2 + K}$, provided $\kappa^2 + K > 0$.  The initial condition will have trivial holonomy only if $L$ is an integral multiple of $2\pi/\sqrt{\kappa^2 + K}$, which does not happen for a typical surface of revolution. An important exception is, of course, the sphere, for which $L = 2\pi/\sqrt{\kappa^2 + K}$. In this case, the holonomy is trivial for every initial condition, reflecting the fact that the sphere admits three independent infinitesimal isometries. If $\kappa^2 + K \le 0$, only multiples of $(1,0,-\kappa)$ have trivial holonomy.

%%%%%%%%%%%%%%%%%%%%%%%%%%%%%%%%%%%%%%%%%%%%%%%%%%%%%%%%%%%%%%%%%%%%%%%%%%%%%%%%%%%

\section{A proof of the Gauss-Bonnet theorem\label{GBT}}

In this section we use Killing transport to give another proof of the Gauss-Bonnet theorem on a surface. Assume that $M$ is a compact, oriented, Riemannian surface. We begin with two lemmas that give some additional properties of Killing transport and infinitesimal isometries along a curve.

\begin{lemma}
\label{EndptChoice} Let $\gamma:[0,L]\to M$ be a $C^2$ curve parameterized by arc length. There exists a non-trivial infinitesimal isometry $(X(t),\xionetwo(t))$ along $\gamma$ such that $X(0)$ is a multiple of $\gamma'(0)$ and $X(L)$ is a multiple of $\gamma'(L)$.
\end{lemma}

\begin{proof}
Define $\mu:T_{\gamma(0)}M\oplus \mathbb{R} \rightarrow T_{\gamma(L)}\oplus \mathbb{R}$ by Killing transport. Let
$V_{t}=\operatorname{span}\{(  \gamma'(  t)  ,0), (0,1)\}\subset T_{\gamma(t)}M \oplus \mathbb{R}$. Then $\mu(V_{0})$
$\cap$ $V_{L}$ is non-trivial by dimension count. Any non-zero element
$(X(L),\xionetwo(L))$ in $\mu(V_{0})$ $\cap$ $V_{L}$ is the terminal
value of an infinitesimal isometry along $\gamma$ with the desired endpoint conditions. 
\end{proof}

\begin{lemma}
\label{ZeroBehavior} Let $\gamma$ be a $C^2$ curve in $M$ parameterized by arc length, and let $(X(t),\xionetwo(t))$ be a non-trivial infinitesimal isometry along $\gamma$. If $X(t_0) = 0$, then the limiting angle between $X(t)$ and $\gamma'(t)$ at $t_0$ is a right angle.
\end{lemma}

\begin{proof}
Since $(X(t),\xionetwo(t))$ is non-trivial, then $\xionetwo(t_0) \ne 0$. Let $\{T, N\}$ be the oriented frame along $\gamma$ with $T = \gamma'$. From (\ref{GenJacEqtn2D}) we have $\nabla_T X = -\xionetwo N$ along $\gamma$, and so the zero for $X(t)$ is isolated. Define the $C^1$ function $r(t)$ and the $C^1$ unit vector field $u(t)$ by
$X(t) = r(t)u(t)$ along $\gamma$ and $r(t) = ||X(t)|| > 0$ for $t_0 < t < t_0 + \epsilon$. Differentiating and setting $t=t_0$ yields $-\xionetwo(t_0)N(t_0) = r'(t_0)u(t_0)$, and the result follows. 
\end{proof}

%*****************************************************************

Suppose $(X(t), \xionetwo(t))$ is an infinitesimal isometry along $\gamma$, where $\gamma$ is a $C^2$ curve parameterized by arc length. Let $e_1$, $e_2$ be an oriented, orthonormal frame along $\gamma$, and let $\omega^1$, $\omega^2$ be the dual frame.  We can write
$$
\gamma' = \cos\tau\, e_{1} + \sin\tau\, e_{2} \qand 
X = \xi^1 e_{1} + \xi^2 e_{2} = r(\cos\theta\,e_1 + \sin\theta\,e_2),
$$
where $\tau$, $r$, and $\theta$ are $C^1$ functions. By Lemma~\ref{ZeroBehavior}, $r$ changes sign at the zeros of $X$. 

Assume that the infinitesimal isometry along $\gamma$ is chosen with endpoint data as in Lemma~\ref{EndptChoice}. Then Lemma~\ref{ZeroBehavior} implies that
$\theta(L^{-})-\theta(0^{+}) = \tau(L^{-}) - \tau(0^{+}) + \frac{\pi}{2}m,$ where $m$ is an integer determined by the transport. (An odd multiple of $\pi/2$ is obtained if and only if $X$ vanishes at exactly one endpoint of $\gamma$.) Note that $m\pi/2$ is the net change of angle between $\gamma'$ and the vector $X/r = \cos\theta\,e_1 + \sin\theta\,e_2$.

The system (\ref{SysOnCurve2D}) implies%
\begin{equation}
\begin{aligned}
d\theta
&= \omega_{2}^{1}-\xionetwo\,r^{-1}(\cos\theta\,\omega^{1}+\sin\theta\,\omega^{2}) \\
&= \omega_{2}^{1}-\xionetwo\,r^{-2}(\xi^{1}\omega^{1}+\xi^{2}\omega^{2}) 
 = \omega_{2}^{1}-\xionetwo\,r^{-2}\left<X,\cdot\right>,
\end{aligned}
\label{LocGaussBonnet}%
\end{equation}
where $\left<X,\cdot\right>$ denotes the 1-form dual to $X$.

Integrate
(\ref{LocGaussBonnet}) along $\gamma$ to obtain%
\begin{equation}%
\begin{aligned}
\textstyle 
\int_{\gamma}\omega_{2}^{1}
- \int_{\gamma}\xionetwo\,r^{-2}\left<X,\gamma'\right>dt
& =\theta(L^{-})-\theta(0^{+}) \\
& \textstyle =\tau(L^{-})-\tau(0^{+})+\frac{\pi}{2}m.%
\end{aligned}
\label{EdgeFormula}%
\end{equation}
Note that the integrand $\xionetwo\,r^{-2}\left<X,\gamma'\right>$ is continuous at the zeros of $X$ because of (\ref{LocGaussBonnet}).

This leads to another proof of the Gauss-Bonnet theorem.

\begin{theorem}
$\int_{M^{2}}K\,dA=2\pi\chi\left(  M^{2}\right)  $ for a compact manifold
$(M,g).$
\end{theorem}

\begin{proof}
Let $\left\{  \triangle^{i}:1\leq i\leq f\text{\/}\right\}  $ be a
triangulation of $(M,g).$ We do \textit{not} assume our
triangulation to be geodesic. Choose an oriented, orthonormal frame on each triangle. Denote the connection form on $\triangle^{i}$ by $(\omega^i)_{2}^{1}$, and recall from Section~\ref{connection} that $\xionetwo$ does not depend on the frame. Write $\partial\triangle^{i}=:\gamma_{1}%
^{i}+\gamma_{2}^{i}+\gamma_{3}^{i}$, where directed edge $\gamma_{\alpha}^{i}$ joins
vertex $p_{\alpha}^{i}$ to vertex $p_{\alpha+1}^{i}$ for $\alpha$ $=$ $1,$ $2,$
$3$, giving $\partial\triangle^{i}$ the positive orientation. (Here, and in the rest of the proof, the index associated with the edge or vertex of a triangle is read 
$\operatorname{mod}3$.) 

Along each edge $\gamma_{\alpha}^{i}$, assign a non-trivial infinitesimal isometry with initial and terminal values as in Lemma~\ref{EndptChoice}.
For adjacent triangles $\triangle^{i}$ and $\triangle^{j}$ having common edge
$\gamma_{\alpha}^{i}=-\gamma_{\beta}^{j}$ (with opposite orientation), assign the same infinitesimal isometry. On edge $\gamma_{\alpha}^{i}$ of triangle $\triangle^{i}$,
(\ref{EdgeFormula}) reads%
\begin{equation}\label{EdgeFormula2}
\textstyle
\int_{\gamma_{\alpha}^{i}}(\omega^i)_{2}^{1}
= \big[\tau_{\alpha}^{i}(L^{-})-\tau_{\alpha}^{i}(0^{+})\big] 
+ \int_{\gamma_{\alpha}^{i}}\xionetwo\,r^{-2}\left<X,\gamma'\right>dt + \frac{\pi}{2}m_{\alpha}^{i}.
\end{equation}
Note that when $\gamma_{\alpha}^{i}=-\gamma_{\beta}^{j}$, we have $m_{\alpha}^{i}=-m_{\beta}^{j}$.

Sum both sides of (\ref{EdgeFormula2}) over all the edges of all triangles. For the common edge of adjacent triangles, the second and the third terms on the right side cancel due opposite orientations.  We have%
$$
\begin{aligned}
&\textstyle
\sum_{i=1}^{f}\int_{\partial\triangle^{i}}(\omega^i)_{2}^{1} 
=\textstyle\sum_{i=1}^{f}\sum_{\alpha=1}^{3}\big[\tau_{\alpha}^{i}(L^{-})-\tau_{\alpha}^{i}(0^{+})\big] \\
&\qquad= \textstyle\sum_{i=1}^{f}\big(2\pi-(\varepsilon_{1}^{i}+\varepsilon_{2}^{i}+\varepsilon_{3}^{i})\big)
= \sum_{i=1}^{f}\big((  \iota^i_{1}+\iota^i_{2}+\iota^i_{3})  -\pi\big),
\end{aligned}
$$
where $\varepsilon_{\alpha}^{i}$ is the exterior angle between $(\gamma_{\alpha-1}^{i})'(L)$ and $(\gamma_{\alpha}^{i})'(0)$ at $p_{\alpha}^{i}$ and $\iota^i_{\alpha}$
$:=$ $\pi$ $-$ $\varepsilon_{\alpha}^i$ is the interior angle of $\triangle^i$ at
$p_{\alpha}^{i}.$ Apply Stokes' theorem to the left side and sum up the right side
to obtain $\int_{M^{2}}K\,dA=2\pi v-\pi f,$ where $v$ is the number
of vertices in the triangulation. Since $e,$ the number of edges, satisfies
$3f=2e,$ we have $2\pi v$ $-\pi f$ $=$ $2\pi(v-e+f) = 2\pi\chi(M)$, and the theorem follows.
\end{proof}

%%%%%%%%%%%%%%%%%%%%%%%%%%%%%%%%%%%%%%%%%%%%%%

\section{Existence of infinitesimal isometries on surfaces\label{existence}}

In this section we determine criteria for the existence of local infinitesimal isometries in dimension two, that is, for local solvability of the system in Proposition \ref{MainSystem}. Our goal is Theorem~\ref{class of inf isom}. While this is a classical result (\cite[Livre VII, Ch.~II]{Darboux} and \cite[pp.~318--322]{Cartan}), we would like to thank R.~Bryant for the discussions with the second author that led to the proof we give here. (Darboux \cite{Darboux} and Cartan \cite{Cartan} address the question of when two surfaces are locally isometric. Once this is determined, they investigate the size (number of parameters) of the family of all local isometries.) Another recent treatment is given in \cite{Atkins}.   A related result on one-parameter families of infinitesimal isometries is given in \cite[pg.~49]{Yano}.

Let $e_1$, $e_2$ be a local orthonormal frame on $M^2$. Consider the Euclidean space $\mathbb{R}^{3}$ of variables
$(\xi^{1},\xi^{2},\xionetwo)$. Then the submanifold of the first jet space
of $\xi$ defined by $\xi_{1}^{1}=\xi_{2}^{2}=0$, $\xionetwo+\xi_{1}^{2}=0$
may be identified with $\calS:=M^{2}\times\mathbb{R}^{3}$.

On $M\times\mathbb{R}^{3}$ consider the Pfaffian system $\theta=(\theta
^{1},\theta^{2},\theta^{3})$ given by
\begin{equation}%
\begin{array}
[c]{lcl}%
\rule[-2mm]{0mm}{6mm}\theta^{1} & = & d\xi^{1}+\xi^{2}\omega_{2}^{1}-\xi
_{2}^{1}\omega^{2},\\
\rule[-2mm]{0mm}{6mm}\theta^{2} & = & d\xi^{2}-\xi^{1}\omega_{2}^{1}+\xi
_{2}^{1}\omega^{1},\\
\rule[-2mm]{0mm}{6mm}\theta^{3} & = & d\xionetwo-K\xi^{2}\omega^{1}+K\xi
^{1}\omega^{2}.
\end{array}
\label{thetas}%
\end{equation}
We check the Frobenius integrability conditions for (\ref{thetas}): By
(\ref{connection forms}) and (\ref{curvature tensor}) we have
\[
d\theta^{1},\,d\theta^{2}\equiv0~~\mod~\theta
\]
and
\[
d\theta^{3}\equiv(K_{1}\xi^{1}+K_{2}\xi^{2})\,\omega^{1}\wedge\omega
^{2}~~~~\mod~\theta,
\]
where $K_{i}=dK(e_{i})$,~~$i=1,2$, so that $dK=K_{1}\omega^{1}+K_{2}\omega
^{2}$. Note that $d\theta^{3}$ mod\,$\theta$ is essentially the curvature of the connection
$\tilde\nabla$ given in equation (\ref{curvature}). We define
$$
T = K_{1}\xi^{1}+K_{2}\xi^{2}
$$
on $\calS$ and consider several cases.

{\it Case 1:} $T$ is identically zero on $\calS$. We have that $T \equiv 0$ on $\calS$ if and only if $K_1$ and $K_2$ vanish identically, that is, $K$ is constant. In this case (\ref{thetas}) is integrable and there exists a three-parameter family of solutions by the Frobenius theorem.

{\it Case 2:} $T$ is not identically zero on $\calS$. In this case we assume $dT \ne 0$ so that $\calS' := \{ T = 0 \}$ is a submanifold of dimension 4. If (\ref{thetas}) has an integral manifold, it will be contained in $\calS'$. Differentiating
$dK=K_{1}\omega^{1}+K_{2}\omega^{2}$, we see by (\ref{connection forms}) that
\begin{equation}%
\begin{array}
[c]{lcl}%
0 & = & d^{2}K\\
& = & (dK_{1}+K_{2}\omega_{2}^{1})\wedge\omega^{1}+(dK_{2}-K_{1}\omega_{2}%
^{1})\wedge\omega^{2}.
\end{array}
\label{ddk}%
\end{equation}
We define $K_{ij}$ so that
\begin{align}
dK_{1}  &  =-K_{2}\omega_{2}^{1}+K_{11}\omega^{1}+K_{12}\omega^{2}%
,\label{dk1}\\
dK_{2}  &  =K_{1}\omega_{2}^{1}+K_{21}\omega^{1}+K_{22}\omega^{2}. \label{dk2}%
\end{align}
By substituting (\ref{dk1}) and (\ref{dk2}) into (\ref{ddk}), we have
$K_{12}=K_{21}$.

Then we have by (\ref{thetas}), (\ref{dk1}) and
(\ref{dk2})
\begin{align*}
dT  &  =\xi^{1}dK_{1}+K_{1}d\xi^{1}+\xi^{2}dK_{2}+K_{2}d\xi^{2}\\
&  \equiv{\ }(K_{11}\xi^{1}+K_{12}\xi^{2}-K_{2}\xionetwo)\,\omega^{1} 
+ (K_{12}\xi^{1}+K_{22}\xi^{2}+K_{1}\xionetwo)\,\omega^{2} \mod~\theta.
\end{align*}
We set
\begin{equation}\label{Ti def}
\left\{
\begin{array}
[c]{lcl}%
T_{1} & = & K_{11}\xi^{1}+K_{12}\xi^{2}-K_{2}\xionetwo,\\
T_{2} & = & K_{12}\xi^{1}+K_{22}\xi^{2}+K_{1}\xionetwo.
\end{array}
\right.
\end{equation}

{\it Sub-case 2.1.} Next we show that $T_1$ and $T_2$ cannot both vanish identically on $\calS'$. Assume, on the contrary, that $T_{1},T_{2}\equiv0$ on $\calS'$. Let $i:\calS' \hookrightarrow \calS$ be the inclusion map. Then $i^{\ast}\theta= (i^{\ast}\theta^{1},i^{\ast}\theta^{2},i^{\ast}\theta^{3})$ has rank
two by Lemma \ref{MainLemma}. Then $\calS'$ is foliated by
two-dimensional integral manifolds and therefore there is a two-parameter
family of solutions. But this is impossible for the following reason. Consider the subset $N := \{\xi^1 = \xi^2 = 0\} \subset \calS'$. Since $T_1$ and $T_2$ vanish on $N$, (\ref{Ti def}) implies that $K_1$ and $K_2$ also vanish on $N$. Then $dT = 0$ on $N$, which contradicts the assumption that $dT \ne0$.

Now we consider the subset $\calS'' := \{T = T_1 = T_2 = 0\}$.  If (\ref{thetas}) has an integral manifold, it will be contained in $\calS''$. Let 
$A=\left(
\begin{smallmatrix}
K_{1} & K_{2} & 0\\
K_{11} & K_{12} & -K_{2}\\
K_{12} & K_{22} & K_{1}%
\end{smallmatrix}
\right)$. 

{\it Sub-case 2.2:} $\det A\ne0$. In this case $\calS''$ is a 2-dimensional submanifold given by $\xi^1 = \xi^2 = \xionetwo = 0$, which means there are no non-trivial infinitesimal isometries.  

{\it Sub-case 2.3:} $A$ has constant rank two. In this case $\calS''$ is a 3-dimensional submanifold of $\calS$. If we have $dT_{1},dT_{2}\equiv0\mod\theta$ at all points of $\calS''$, then Lemma \ref{MainLemma} and the Frobenius theorem imply that $\calS''$ is foliated by two-dimensional integral manifolds and therefore there exists a one-parameter family of solutions. Thus, we seek to express $dT_{1},dT_{2}\mod\theta$ on $\calS''$ in terms of $\omega^1$, $\omega^2$, $\omega^1_2$. 
Differentiating (\ref{dk1}) and (\ref{dk2}) we have
\begin{equation}\label{ddk1}
\begin{aligned}
0  &  =d^{2}K_{1}\\
&  =(dK_{11}+2K_{12}\,\omega_{2}^{1})\wedge\omega^{1} 
+ (dK_{12}+K_{22}\,\omega_{2}^{1}-K_{11}\,\omega_{2}^{1})\wedge\omega^{2} \\
&\qquad- K_{2}K\,\omega^1\wedge\omega^{2}
\end{aligned}
\end{equation}
and
\begin{equation}\label{ddk2}
\begin{aligned}
0  &  =d^{2}K_{2}\\
&  =(dK_{12}+K_{22}\,\omega_{2}^{1}-K_{11}\,\omega_{2}^{1})\wedge\omega^{1} 
      +(dK_{22}-2K_{12}\,\omega_{2}^{1})\wedge\omega^{2} \\
     &\qquad + K_{1}K\,\omega^{1}\wedge\omega^{2}.
\end{aligned}
\end{equation}
We define $K_{ijk}$ so that
\begin{align}
dK_{11}  &  =-2K_{12}\,\omega_{2}^{1}+K_{111}\,\omega^{1}+K_{112}\,\omega^{2},\label{dk11}\\ 
dK_{12}  &  =(K_{11}-K_{22})\,\omega_{2}^{1}+K_{121}\,\omega^{1}+K_{122}\,\omega^{2}, \label{dk12} \\
dK_{22}  & =2K_{12}\,\omega_{2}^{1}+K_{221}\,\omega^{1}+K_{222}\,\omega^{2}. \label{dk22}
\end{align}
Substituting (\ref{dk11}), (\ref{dk12}), and (\ref{dk22}) into (\ref{ddk1}) and (\ref{ddk2}) we have 
\begin{equation}\label{Kijk symmetries}
K_{112}=K_{121}-K_{2}K \qquad\text{and}\qquad K_{122} = K_{221} + K_1 K.
\end{equation}
We have by (\ref{thetas}), (\ref{dk11}), (\ref{dk12}), (\ref{dk22}), and (\ref{Kijk symmetries})
\begin{align*}
dT_{1} & \equiv (K_{111}\xi^{1}+K_{112}\xi^{2}-2K_{12}\xi_{2}^{1})\,\omega^{1} \\
&\qquad + \big(K_{121}\xi^{1}+K_{122}\xi^{2}+(K_{11}-K_{22})\xionetwo\big)\,\omega^{2}~~\mod~\theta
\end{align*}
and
\begin{align*}
dT_{2} & \equiv \big(K_{121}\xi^{1}+K_{122}\xi^{2}+(K_{11}-K_{22})\xi_{2}^{1}\big)\,\omega^{1}\\
&  \qquad+(K_{221}\xi^{1}+K_{222}\xi^{2}+2K_{12}\xionetwo)\,\omega^{2}~~\mod~\theta.
\end{align*}
Note that the terms containing $\omega^1_2$ drop out because $T_1 = T_2 = 0$ on $\calS''$.

We summarize the discussions of this section in the following theorem.

\begin{theorem}
\label{class of inf isom}
Let $M$ be a Riemannian manifold of dimension $2$.  Assume $e_1$, $e_2$ is an orthonormal frame on $M^2$. Relative to this frame, define 
$$\mathbf{K}=
\begin{pmatrix}
K_{1} & K_{2} & 0\\
K_{11} & K_{12} & -K_{2}\\
K_{12} & K_{22} & K_{1}\\
K_{111} & K_{112} & -2K_{12}\\
K_{121} & K_{122} & K_{11}-K_{22}\\
K_{221} & K_{222} & 2K_{12}%
\end{pmatrix}
.
$$

\begin{itemize}
\item[(i)] If $\mathbf{K}$ has rank 0, there exists a three-parameter family
of infinitesimal isometries,

\item[(ii)] $\mathbf{K}$ cannot have rank 1, and there does not exist a two-parameter family
of infinitesimal isometries,

\item[(iii)] If $\mathbf{K}$ has rank 2 and $(K_{1},K_{2})\ne0$, there exists
a one-parameter family of infinitesimal isometries,

\item[(iv)] If $\mathbf{K}$ has rank 3, there exists only the trivial
infinitesimal isometry.
\end{itemize}
\end{theorem}

%%%%%%%%%%%%%%%%%%%%%%%%%%%%%%%%%%%%%%%%%%%%%%%%%%%%%%%%%%

\end{document}